\documentclass[12pt]{article}

\usepackage{amsmath,amssymb,amsfonts,amsthm}
\usepackage{eucal}

\usepackage{verbatim} 

\usepackage{graphicx,psfrag} 

\usepackage{color}

\usepackage{ifthen} 

\newcommand{\baseenvskip}{\baselineskip 6.53mm}

\newtheorem{thm}{Theorem}[section]

\newtheorem{lemma}[thm]{Lemma}

\newtheorem{conjecture}[thm]{Conjecture}

\newtheorem{corollary}[thm]{Corollary}

\theoremstyle{remark}
\newtheorem{rmk}{Remark}[section]
\newenvironment{remark}{\begin{rmk}\rm\baseenvskip}{\end{rmk}}

\numberwithin{equation}{section} \numberwithin{thm}{section}
\numberwithin{rmk}{section} \numberwithin{figure}{section}

\newcommand{\jl}{$\frac{}{}$}

\newcommand{\ds}{\displaystyle}
\newcommand{\QED}{\hfill $\Box$}

\newcommand{\set}[1]{\left\{#1\right\}}

\newcommand{\n}{\nabla}
\newcommand{\Hv}{\mathcal{H}_v}

\newcommand{\Ov}{\mathcal{O}_v}
\newcommand{\R}{\mathbb{R}}
\newcommand{\Rn}{\mathbb{R}^n}
\newcommand{\C}{\mathbb{C}}
\newcommand{\Cn}{\mathbb{C}^n}
\newcommand{\Z}{\mathbb{Z}}
\newcommand{\p}{\partial}

\newcommand{\gv}{\gamma_v}

\newcommand{\eps}{\varepsilon}

\newcommand{\noi}{\noindent}
\newcommand{\finv}{f^{-1}}

\newcommand{\Gn}{\Gamma^{n-1}}
\newcommand{\W}{W^{n-1}}
\newcommand{\We}{W^{n-1}(\varepsilon)}
\newcommand{\Wt}{W^{n-1}(t)}

\newcommand{\Ze}{Z^{n-1}(\varepsilon)}
\newcommand{\Zt}{Z^{n-1}(t)}
\newcommand{\Zr}{Z^{n-1}(R)}
\newcommand{\Zn}{Z^{n-1}}
\newcommand{\el}[1]{e(#1)_{\ell}}

\newcommand{\vf}{\textrm{Lk}}

\newcommand{\ga}{\gamma}

\newcommand{\Nv}{\mathcal{N}_v}

\newcommand{\U}{\mathcal{U}}
\newcommand{\vez}{u_0}
\newcommand{\veu}{u_1}

\newcommand{\lk}{\textsf{Lk}}

\newcommand{\pic}{1}

\begin{document}

\title{Foliations and Global Inversion}


\author{E. Cabral Balreira\\{\small
Department of Mathematics}\\{\small Trinity University}\\ {\small San Antonio, TX 78212}\\{\small ebalreir@trinity.edu}}


\maketitle

\begin{abstract}
We consider topological conditions under which a locally invertible
map admits a global inverse. Our main theorem states that a local
diffeomorphism $f: M \to\mathbb{R}^n$ is bijective if and only if
$H_{n-1}(M)=0$ and the pre-image of every affine hyperplane is
non-empty and acyclic. The proof is based on some geometric
constructions involving foliations and tools from intersection
theory. This topological result generalizes in finite dimensions the
classical analytic theorem of Hadamard-Plastock, including its
recent improvement by Nollet-Xavier. The main theorem also relates to a conjecture of the aforementioned authors, involving the well known Jacobian Conjecture in algebraic geometry.
\end{abstract}

\section{Introduction}\label{chpt:Intro}

In this paper we are concerned with the problem of finding
topological conditions ensuring that a local diffeomorphism is
bijective. A classical result in this direction is the well-known
Hadamard-Plastock Theorem (see~\cite{Nirenberg} and
\cite{Plastock}). It states that a Banach space local diffeomorphism
$f:X\to X$ is bijective provided
 \begin{equation}\label{eq:Had}
\inf_{x\in X} \|Df(x)^{-1}\|^{-1} > 0.
\end{equation}

The proof of the Hadamard-Plastock theorem follows from simple
arguments involving covering spaces. In recent years new topological
and geometric ideas have been introduced in the subject of global
invertibility, pushing the field in different directions (see, for
instance, \cite{NolletXavierTaylor}, \cite{NolletXavier3}, \cite{NolletXavier1}, \cite{NolletXavier2}, \cite{Rabier2}, \cite{SmythXavier}, \cite{Xavieridentity}, and \cite{Xavier}). The emerging picture reveals that global invertibility is also influenced by more subtle topological phenomena. In \cite{NolletXavier1}, Nollet and Xavier established a substantial improvement to the Hadamard-Plastock theorem when $\dim X<\infty$. Using degree theory, they showed in \cite{NolletXavier1} that a local diffeomorphism $f:\Rn\to\Rn$ is bijective if there exists a complete Riemannian metric $g$ on $\Rn$ such that,
 \begin{equation}\label{eq:NX1}
   \forall~v\in S^{n-1},~ \inf_{x\in \mathbb{R}^n} \|Df(x)^* v\|_g >0.
\end{equation}

Notice that \eqref{eq:NX1} is an improvement over \eqref{eq:Had}
since
 \[
\|Df(x)^{-1}\|^{-1} = \|Df(x)^{\ast-1}\|^{-1} =
\ds\inf_{|v|=1}\|Df(x)^{\ast}v\|.
\]

\noi Furthermore, it is easy to produce examples that satisfy
\eqref{eq:NX1} but not \eqref{eq:Had}. Arguments from elementary
Morse theory (see \cite[p.112]{Palmeira}) show that if \eqref{eq:NX1}
holds, then the pre-images of affine hyperplanes $H$ must satisfy
{\nolinebreak
$\finv(H)\times\nolinebreak\R\nolinebreak\cong\nolinebreak\Rn$}
(note that $Df(x)^{\ast}v = \n \langle f(x),v\rangle$). In
particular, by the K\"{u}nneth formula, $\finv(H)$ is
\textit{acyclic} (recall that a topological space is called acyclic
if it has the homology of a point).

In this paper we show that the above mentioned analytical results
are but a manifestation of a topological phenomenon.


\begin{thm}\label{thm:characterization result}
A local diffeomorphism $f:\mathbb{R}^n\to\mathbb{R}^n$ is bijective
if and only if the pre-image of every affine hyperplane is non-empty
and acyclic.
\end{thm}


In Section~\ref{sec:Proof of INJ} we will point out a connection between the above theorem and the Jacobian Conjecture in algebraic geometry. The non-trivial half of Theorem~\ref{thm:characterization result} consists in establishing injectivity and surjectivity. Its proof is based on some geometric constructions involving foliations, and the computation of intersection numbers of certain chain complexes. Theorem~\ref{thm:characterization result} also allows for an analytic corollary that is stronger than the results in \cite{NolletXavier1}, in the sense that one can choose the metric to suit the unit vector $v$.

\nocite{Xavierrigidity}\nocite{Hadamard}


\section{Preliminaries}\label{chpt:Combinatorial Topology}
Given a compact smooth manifold $M^n$ and a finite cover, we would
like to have a systematic way to describe the intersections of the
sets in the cover. Likewise, once a point is given we want to
describe exactly all the sets in the cover that contain the given
point. To this end, we will consider a triangulation of $M$ and view
the top dimensional cells as the sets of the covering. We set our
notation as follows. Denote by $T(M)$ a triangulation on $M$ (whose
existence is guaranteed by \cite{Whitehead}) and let $e(k)_j$ be the
$j^{th}$ $k$-cell of $T(M)$. The set of indexes of $k$-cells will be
denoted by $E(k)\subset \mathbb{N}$. Also, given a triangulation
$T(M)$, let $T_k(M)$ be the $k$-\textit{skeleton} of $M$. Whenever
the context is clear, we will refer to the triangulated space simply
as~$M$.

This combinatorial approach allow us to easily address the
properties we mentioned above. For instance, given a simplex
$e(k)_j$, the \textit{star} of $e(k)_j$ describes all the simplexes
that contain $e(k)_j$. In our results, we will be interested in
finding all the $(k+1)$-simplexes that contain $e(k)_j$. This is
easily accomplished by looking at the vertices of the \textit{link}
of $e(k)_j$, denoted by $\lk\left(e(k)_j\right)$.

We now review the basic definitions from intersection theory. We
define in $M^n$ the intersection number (mod 2) between $A^p$ a
$p$-cycle and $B^{q}$ a $q$-cycle, where $p+q=n$ by $\#(A^p,B^{q})$.
We note that when $A^p$, $B^{q}$ represent transverse submanifolds,
then $\#(A^k,B^{n-k})$ represents the number of geometric
intersections mod~2. The property that we highlight is that
intersection number depends only on the homology class. For details
and formal definitions we refer the reader to \cite{Schwartz}.

Finally, we can also define linking numbers between cycles. Let
$X^p$ and $Y^{q-1}$ be two nonintersecting cycles in $\Rn$ with
$p+q=n$. For $Z^{p+1}$ a bounding chain of $X^p$, i.e. $\p Z^{p+1} =
X^p$, we define the \emph{linking number} between $X^p$ and
$Y^{q-1}$ as
 \begin{equation}\label{eq: def of linking number}
\vf(X^p,Y^{q-1}) = \#(Z^{p+1}, Y^{q-1}),
\end{equation}

\noi which is independent of the choice of the bounding chain of
$X^{p}$.


\section{Injectivity}\label{sec:Proof of INJ}
Let us consider a local diffeomorphism $f:M\to\Rn$, where $M$ is a
smooth connected manifold. Our goal is to understand under which
topological conditions the map $f$ is injective. There is a conceptual link between injectivity and connectedness. For instance, it is clear that a locally invertible map is injective if and only if the pre-image of every 0-dimensional affine subspace (i.e., a point) is connected (possibly empty). An analogous statement can be made if one goes one dimension higher and considers lines instead of points, that is, a locally invertible map is injective if the pre-image of every line is connected.

In view of these observations, Nollet and Xavier \cite{NolletXavier1} made the following conjecture.

\begin{conjecture} A local diffeomorphism $f: \mathbb{R}^n \rightarrow \mathbb{R}^n$ is injective if the pre-images of every affine hyperplane is connected (possibly empty).
\end{conjecture}

At the present time this conjecture remains open and its significance is better seen in Algebraic Geometry where it would provide a positive answer for the Jacobian Conjecture (recall that the Jacobian Conjecture states that a polynomial local biholomorphism $F:\C^n\to\C^n$ is invertible, see \cite{BassCW}, \cite{Essen}). Indeed, if $F:\Cn\to\Cn$ is a polynomial local biholomorphism, and $H\subset \Cn$ is a real hypersurface foliated by complex hyperplanes $V$, then by a Bertini type theorem $F^{-1}(V)$ is connected for a generic $V$ (see \cite{Schinzel}, Cor.~1 of Theorem~3.7). From this one can easily check that $F^{-1}(H)$ is connected and hence one would establish the Jacobian Conjecture.

The result below establishes a weaker version of the Nollet-Xavier conjecture, where connectedness is replaced by acyclicity.


\begin{thm}\label{thm:inj_result}
A local diffeomorphism $f:\mathbb{R}^n\to\mathbb{R}^n$ is injective
if the pre-image of every affine hyperplane is either empty or
acyclic.
\end{thm}

In fact, we observe that we may weaken the hypotheses of Theorem~\ref{thm:inj_result} to obtain the following stronger result. We say that an affine hyperplane $H\subset \Rn$ is parallel to a line $\ell$ in $\Rn$ provided that $\ell\cap H=\emptyset$ or $\ell\subset H$.

\begin{thm}\label{thm:injetividade general result for a manifold}
For $n\geq 3$, let $f:M\to\Rn$ be a local diffeomorphism where $M$
is a (necessarily non-compact) connected manifold with $H_{n-1}(M)=
0$. If there exists a line $\ell$ in $\Rn$ such that the pre-image
of every affine hyperplane parallel to $\ell$ is either empty or
acyclic, then $f$ is injective.
\end{thm}

The proof of Theorem~\ref{thm:injetividade general result for a
manifold} is based on geometric constructions of chain complexes,
the computation of the intersection number between these objects,
and the maximal lift of lines. Since the computation of intersection
numbers is done with objects belonging to the domain of~$f$, we need
to require the extra assumption on the homology of $M$. Observe that
the cases in Theorem~\ref{thm:injetividade general result for a
manifold} when $n=1,2$ are trivially true without any extra
assumptions on $M$.

We stress that in our arguments we will only require the existence
of \emph{local lifts}. In fact, by Hadamard-Plastock
Theorem~\cite{Plastock}, if all lines admit \emph{global lifts} the
map is already bijective. We refer to a local lift of a line
$\ell=\set{tw|w\in\Rn, t\in\R}$ with respect to $f$ as a path
$\alpha:(-\eps,\eps)\to M$, $\eps>0$ such that
$f\left(\alpha(t)\right)=tw$. Observe that by the Inverse function
theorem, if $f$ is a diffeomorphism a local lift of a line always
exists in the above sense. We say that $\alpha:(-\delta,\delta)\to
M$ is the \emph{maximal lift} of $\ell$ if
$\delta=\sup\{\eps|\alpha\textrm{ admits a local lift for
}\eps>0\}$. Furthermore, $\ell$ has a \emph{global lift} if its
maximal lift satisfies $\delta=\infty$. Finally, what is important
for us is the fact that the maximal lift of a line is properly
embedded in the domain. In our notation, $\alpha$ is properly
embedded if it leaves every compact set of $M$ as $|t|$ increases to
$\delta$.

\subsection{Beginning of the proof of Theorem~\ref{thm:injetividade general result for a manifold}}

Assume that there is a point $p$ in the image of $f$ with at least
two distinct points $q_0$ and $q_1$ in its pre-image. Since
translations do not change any of the hypotheses, we assume for
simplicity that $p=0$. Our goal is to construct a $(n-1)$-cycle
$\Gn$ so that the intersection number of $\Gn$ with the maximal lift
of the line $\ell$ passing through the origin will necessarily be
zero, as the maximal lift is properly embedded. A simple argument
will then show that $f$ must have a critical point along the lift,
thus establishing the desired contradiction.

First, we give an outline of the proof. Consider $\eps>0$ so that
the ball $V=B(0;\eps)\subset f(M)$ has diffeomorphic pre-images
$U_0$ and $U_1$ around $q_0$ and $q_1$, respectively. Next, let $Y$
be the $(n-1)$-\textit{equatorial disk} of $V$ determined by $\ell$,
that is, the intersection of the orthogonal hyperplane to $\ell$ and
$V$. The cycle $\Gn$ we seek will be constructed to resemble a
topological cylinder that connects the induced equatorial disks of
$U_0$ to $U_1$, denoted by $X_0$ and $X_1$, respectively. We
construct~$\Gn$ as follows. Take a hyperplane parallel to $\ell$
which intersects $\p V$ tangentially at $v$ and is denoted by $\Hv$.
As we change the hyperplane $\Hv$ by moving it around $\p V$, the
pre-images $u_i\in \p U_i$ of $v$, for $i=0,1$, can be continuously
connected by paths in $\finv(\Hv)$, at least for nearby hyperplanes.
In this way we construct small \emph{lateral} pieces of $\Gn$.

One then tries to put together all those local data. In so doing,
one is forced to consider the situation where, for a fixed
hyperplane $\Hv$, there are multiply-defined paths joining the same
pre-images of points in $\finv(\Hv)$. Whenever this occurs, the
topological hypotheses that $\finv(\Hv)$ is acyclic will be used to
\emph{fill in the gaps}. See Fig.~\ref{fig:sausage by evolution} for
a depiction of this process when $n=3$.

\begin{figure}[h]
\begin{center}
\ifthenelse{\equal{\pic}{1}}{
\begin{psfrags} \psfrag{p}{$p$}
\psfrag{q0}{$q_0$} \psfrag{q1}{$q_1$} \psfrag{u1}{$\veu$}
\psfrag{u0}{$\vez$} \psfrag{pX0}{$\p X_0$} \psfrag{pX1}{$\p X_1$}
\psfrag{H}{$\Hv$} \psfrag{f-1H}{$f^{-1}(\Hv)$} \psfrag{ell}{$\ell$}
\psfrag{Wn1}{$W^{n-1}$} \psfrag{V}{$V$} \psfrag{v}{$v$}
\psfrag{pY}{$\p Y$}
\includegraphics{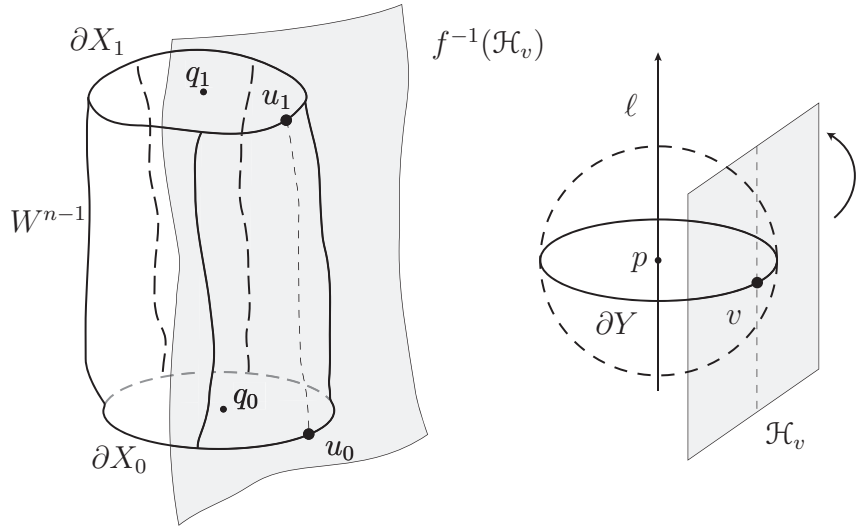}
\end{psfrags}
} {
\[\fbox{\rule[0.75in]{2.5in}{0in}outline of sausage in intro}\]
} \caption{Construction of paths connecting $\p X_0$ to $\p X_1$ by
the revolution of affine hyperplanes.} \label{fig:sausage by
evolution}
\end{center}
\end{figure}

In order to determine how the lateral pieces will fit together and
how such gaps should be filled, we consider a combinatorial
decomposition of $\p Y$ in terms of a triangulation. Here we observe
that $\p Y\cong S^{n-2}$, so such triangulation always exist. The
process of putting together the pieces of $\Gn$ will be done in
steps according to the dimension of the carrier of each point. More
precisely, first we consider a point and the $(n-2)$-cells it may
possibly belong and construct chain complexes that correspond to the
lateral pieces indicated above. Next, points that belong to the
lower dimensional skeleton of $\p Y$ will be consider more than once
in the initial step. Hence, in the following step we consider the
$(n-3)$-skeleton of $\p Y$ and determine the bounding chains
according to the higher dimensional cells that contain it. We repeat
this process until we consider the $0$-skeleton of $\p Y$. The
existence and properties of $\Gn$ are established in the following
lemma.

\begin{lemma}\label{lemma:Chain complex sausage of revolution}
There exists a geometric (singular) chain complex $\Gn\in
S_{n-1}(M)$ that may be represented as $\Gn = X_0 + \W +X_1$ such
that $\W$ is a chain complex with $\p W = \p X_0 + \p X_1$ and for
all $q\in \text{supp\,}\W$ (in its image in $M$), there exists
$v\in\p Y$ so that $q\in \finv(\Hv)$.
\end{lemma}

The proof of Lemma~\ref{lemma:Chain complex sausage of revolution}
follows the outline above where we will construct all the singular
chain complexes of $\W$ and the attaching maps. This argument uses
ideas from combinatorial topology and we postpone it until next
section. We proceed to establish Theorem~\ref{thm:injetividade
general result for a manifold}, but first we remark that we are
interested in the existence of a \textit{geometric} intersection
(i.e., number of points in the set theoretical intersection) between
$\Gn$ and the maximal lift of $\ell$. Therefore we consider
intersection numbers and homology with $\Z_2$ coefficients, thus
avoiding heavier notational concerns regarding orientation and
leaving the proof simpler and more geometric.

Assuming $\Gn$ is constructed as in Lemma~\ref{lemma:Chain complex
sausage of revolution}, we compute the intersection number of $\Gn$
and the maximal lift of $\ell$ starting at $q_0$ which we denote by
$\ga$. We claim that $\ga$ must intersect $\Gn$ in another point
besides $q_0$ and we will show that it is $q_1$. First, we see that
$\ga$ is properly embedded in $M$. Indeed, if we decompose $\gamma$
as $\gamma_-\wedge\gamma_+$ as the maximal lift of $\ell$ in the
negative and positive direction, respectively, starting at $q_0$. It
is then clear that $\gamma_-$ and $\gamma_+$ are not entirely
contained in any compact subset of $M$, otherwise the lift would not
be maximal.

Now, the fact that $\Gn\in S_{n-1}(M)$, i.e., a cycle and
$H_{n-1}(M) =0$ implies that there exists a bounding singular chain
$\Sigma^n$, with $\p \Sigma^n = \Gn$ and a compact set $K$ so that
$\Sigma^n\subset K$. Thus $\Gn$ is a representative of the trivial
element in $H_{n-1}(M,M-K)$. We also have that $\gamma\in
H_1(M,M-K)$ and from the fact that intersection numbers depend only
on the homology class, we have
 \begin{equation}\label{eq:zero cycles}
 \#(\Gn,\gamma) = 0.
\end{equation}

\noi Indeed, $H_{n-1}(M,M-K) = H_{n-1}(M/M-K)$ and since
$\Sigma^n\subset K$, we have $\Gn\sim0$ in $H_{n-1}(M/M-K)$ as well.

From Lemma~\ref{lemma:Chain complex sausage of revolution} and by
definition of intersection numbers, we can write \eqref{eq:zero
cycles} as:
 \begin{equation}\label{eq:intersection of two closed chains mod2}
\begin{array}{rccccc}
0 = & \#(\Gn,\gamma) &&&&\\
= & \#(X_0,\ga) &+& \#(W^{n-1},\ga) &+& \#(X_1,\ga)\\
= & 1 &+& 0 &+& \#(X_1,\ga),
\end{array}
\end{equation}

\begin{figure}[h]
\begin{center}
\begin{psfrags}
\psfrag{p}{$p$} \psfrag{pprime}{$v$} \psfrag{q0}{$q_0$}
\psfrag{q1}{$q_1$} \psfrag{q0prime}{$u_0$} \psfrag{q1prime}{$u_1$}
\psfrag{U0}{$U_0$} \psfrag{U1}{$U_1$} \psfrag{H}{$H$}
\psfrag{f-1H}{$f^{-1}(H)$} \psfrag{ell}{$\ell$}
\psfrag{Wn1}{$\Gamma^{n-1}$} \psfrag{V}{$V$}
\psfrag{line}{$f^{-1}(\ell)$}
\includegraphics{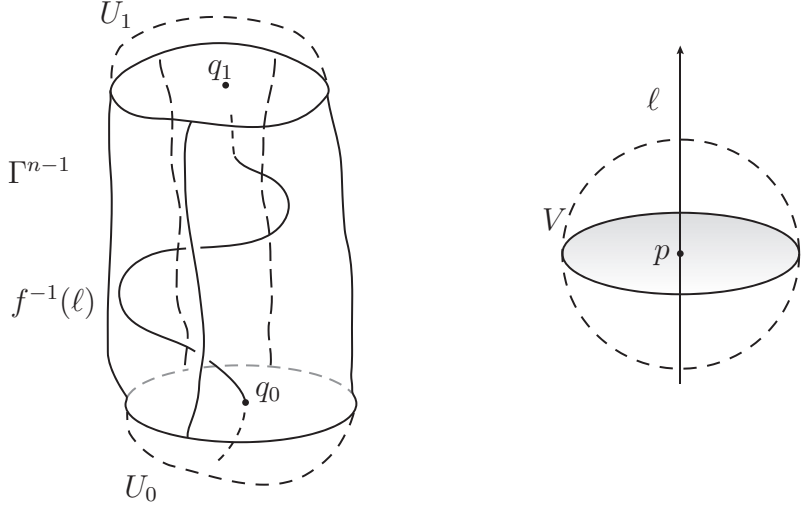}
\end{psfrags}
\caption{Construction of a closed chain complex $\Gn$ by revolving
affine hyperplanes.} \label{fig:Intersection of Chain complexes}
\end{center}
\end{figure}

\noi where the first term is 1 since $f$ is a local diffeomorphism
and the images of $\ga$ and $X_i$ are orthogonal and the second term
is zero since $\ga\cap\W=\emptyset$. Therefore, it must be that
$\#(X_1,\ga)=\nobreak1$. In particular, $\ga\cap X_1\neq\emptyset$
and by the choice of $\eps$, it must be that $\ga\cap
X_1=\set{q_1}$. For a geometric depiction see
Fig.~\ref{fig:Intersection of Chain complexes}.

Finally let $\alpha\subset\finv(\ell)$ be the path segment from
$q_0$ to $q_1$. The image of $\alpha$ is a loop in $\ell$ that has a
point $\widehat{p}\in\ell\cap f(M)$ that is furthest from $p$. Now
it is clear that $f$ fails to be locally invertible at the
corresponding pre-image of $\widehat{p}$, giving us the desired
contradiction. Therefore $f$ must be injective. \QED


\subsection{Reassemblage of Hyperplanes and a Chain Complex
Construction}

We now establish Lemma~\ref{lemma:Chain complex sausage of
revolution}, needed to complete the proof of
Theorem~\ref{thm:injetividade general result for a manifold}. While
outlining the construction of $\Gn$ earlier, we encountered a key
problem which simply put is attributed to the lack of uniqueness on
the choice of the path used to connect the pre-images of a point in
$\p Y$. In our construction this is reflected as follows: although
each path may be defined continuously within a neighborhood of a
fixed point, as we consider the intersection of two neighborhoods
there will possibly be two choices of paths. We claim that whenever
ambiguity occurs, we may use the hypotheses of acyclicity of the
pre-images of hyperplanes to define chain complexes to circumvent
this problem.

We do this by considering a triangulation of $\p Y$ with
sufficiently small mesh to be determined during the proof.
Heuristically, we view a neighborhood of a generic point as the top
dimensional cell containing it and the triangulation will provide a
way to keep track of the intersection of the multiple neighborhoods.
Let $e(k)_j$ be the cells of such triangulation, where
$k=0,\ldots,n-2$ denotes the dimension of each cell and $j\in
E(k)\subset \mathbb{N}$ is the indexing set of the $k$-cells. From
the initial choice of $\eps>0$, we may also define an induced
triangulation via the local diffeomorphism on $\p X_0$ and $\p X_1$
with cells $e(k)_j^0$ and $e(k)_j^1$, respectively.

We construct $\W$ in $n-1$ steps which we enumerate from 0 to $n-2$.
In step $k$ we consider points in the $(n-2-k)$-skeleton of $\p Y$,
denoted by $\p Y_{(n-2-k)}$, and show that the possibly multiply
defined chains are obtained by looking at all the higher dimensional
cells containing such points and that these chains give rise to a
cycle. Then by using the acyclicity hypotheses, we have that such
cycle can be realized as the boundary of another chain complex which
will be the building blocks of $\W$.

\underline{Step $0$:} The initial process is analogous to what has
been outlined before, but to establish our notation we provide the
formal argument. Given the initial triangulation of $\p Y$, take
$v\in \el{n-2}$ and let $u_i$~be the pre-image of $v$ in $\p X_i$
for $i=0,1$. From the connectedness hypotheses of $\finv(\Hv)$,
there is a path $W_{\ell}^1(v)\subset\finv(\Hv)$ joining $u_0$ to
$u_1$, that is, $W_{\ell}^1(v)$ is a 1-chain with $\p W_{\ell}^1(v)
= u_0+u_1$.

Next, we can continuously modify $W_{\ell}^1(v)$ for all points in a
neighborhood of $v\in \p Y$. This follows because $W_{\ell}^1(v)$ is
compact and $f$~is a local diffeomorphism. By repeating this
construction for every point in $\p Y$, we obtain a cover of $\p Y$
from which we extract a finite subcover as $\p Y$ is compact. Then
take finitely many barycentric subdivisions of $\p Y$ until its mesh
is smaller than the minimum diameter of the subcover.

Finally we redo the assignment of $W^1_{\ell}(v)$ for each
$v\in\el{n-2}$ using the newly obtained triangulation. This has the
property that for each $(n-2)$-cell we may define a $(n-1)$-chain
complex denoted by $W^1_{\ell}\times\el{n-2}$ from the continuous
family of paths for each $\ell\in E(n-3)$.

\underline{Step $1$:} In this next step, we consider points in the
$(n-3)$-skeleton of $\p Y$ as these are the points which we possibly
assigned two different 1-chain complexes in the previous step. For
$v\in\el{n-3}$, we may identify all the $(n-2)$-cells that contain
$\el{n-3}$ by looking at the vertices of $\lk(\el{n-3})$. In this
case, we have precisely two points as $\el{n-3}$ belongs to exactly
two top dimensional cells say, $e(n-2)_1$ and $e(n-2)_2$. From the
previous step, we constructed two possibly distinct chain complexes
$W^1_{1}(v)$ and $W^1_{2}(v)$ contained $\finv(\Hv)$ joining $\vez$
to $\veu$. If it is the case they are already the same, we are done.
Otherwise, consider the 1-chain $U^1_{\ell}(v) = W^1_{1}(v) +
W^1_{2}(v)$. We claim $U^1_{\ell}(v)$ is a cycle. Indeed, $\p
U^1_{\ell}(v) = \p W^1_{1}(v) + \p W^1_{2}(v) = \vez + \veu + \vez +
\veu = 0$, since we are using $\Z_2$-coefficients. From the
hypotheses that $\finv(\Hv)$ is acyclic, we have that
$U^1_{\ell}(v)$ is the boundary of a 2-chain denoted by
$W^2_{\ell}(v)$.

Now, using the fact that $f$ is a local diffeomorphism and
$W^2_{\ell}(v)$ is compact, we can continuously define
$W^2_{\ell}(u)$ for all $u$ in a neighborhood of $v$ in $\p
Y_{(n-3)}$. Note that in this step we are only considering points in
the $(n-3)$-skeleton. Therefore we obtain a cover of $\p Y_{(n-3)}$
which by compactness we extract a finite subcover. Next, we iterate
finitely many barycentric subdivisions of the triangulation on $\p
Y$ until its mesh is smaller than the minimum diameter of the
subcover. We then redo the construction of the chain complexes up to
this point in step 0 and 1 using the new triangulation. We do this
so the 2-chain complex defined above can be continuously assigned
for each point within a $(n-3)$-cell and we obtain a $(n-1)$-cell
denoted by $W_{\ell}^{2}\times e(n-3)_{\ell}$ for each $\ell\in
E(n-3)$.

\underline{Step $k$:} For a generic step $k$ ($1< k \leq n-2$), we
consider points in the $(n-2-k)$-skeleton of $\p Y$. For
$v\in\el{n-2-k}$, we look at the $(n-1-k)$-cells that contain
$\el{n-2-k}$. This is the case because in the previous step $k-1$,
we have defined $k$-chains $W^k_{i}(v)$ over $v$ these cells that
$v$ belong, for some $i$. A systematic way to consider these cells
is to look at the vertices of $\lk(\el{n-2-k})$. Let us assume that
those are $e(n-1-k)_1, e(n-1-k)_2, \ldots, e(n-1-k)_j$. We now
define $U^k_{\ell}(v) = W^k_{1}(v) + \cdots + W^k_{j}(v)\subset
\finv(\Hv)$ and we claim that $U^k_{\ell}(v)$ is a cycle. Indeed,
 \begin{equation}\label{eq:computaion in chp3 of closed in step k}
\begin{split}
\p U^k_{\ell}(v) &= \p \left(\sum_{i=1}^j W^k_{i}(v)\right) =
\sum_{i=1}^j  \p W^k_{i}(v)\\
& = \sum_{i=1}^j U^{k-1}_{i}(v) = \sum_{\ell'} W^{k-1}_{\ell'}(v)
\end{split}
\end{equation}

\noi where the chains $U^{k-1}_{i}(v)$ were constructed in the
previous step in a similar manner and $\ell'$ corresponds to the
index of all $(n-k)$-cells that contains $v$.

The chain $U^{k-1}_{i}(v)$ is formed by looking at all the
$(n-k)$-cells that contain $e(n-1-k)_i$ and hence will contain
$\el{n-2-k}$. Therefore, these $(n-k)$-cells can also be determined
by looking at the edges of $\lk(\el{n-2-k})$. Observe that for a
fixed $i$, as we look at the chains of type $W^{k-1}_{\ell'}(v)$
that comprise $U^{k-1}_{i}(v)$ we can alternatively look at the
collection of edges in $\lk(\el{n-2-k})$ that make up
$U^{k-1}_{i}(v)$ and the chains $W^{k-1}_{\ell'}(v)$ will be the
vertices of such edges. However, because each edge contains exactly
two vertices, as we do this for all $i$ each term in the last
summation in~\eqref{eq:computaion in chp3 of closed in step k}
appears twice. Since our computation uses $\Z_2$ coefficients, we
have ~\eqref{eq:computaion in chp3 of closed in step k} is zero
establishing that $U^k_{\ell}(v)$ is a cycle.

Again, from the hypotheses that $\finv(\Hv)$ is acyclic, we find a
bounding $(k+1)$-chain $W^{k+1}_{\ell}(v)\subset \finv(\Hv)$ of
$U^k_{\ell}(v)$, that is, $\p W^{k+1}_{\ell}(v) = U^k_{\ell}(v)$.
Next, an analogous argument as in step 1 is used to find a
neighborhood of $v$ in $\p Y_{(n-2-k)}$ where the assignment of
$W^{k+1}_{\ell}(v)$ is continuous for all points within it. This
follows from the local diffeomorphism of $f$ and compactness if
$W^{k+1}_{\ell}(v)$. This induces a cover of $\p Y_{(n-2-k)}$ and by
compactness we extract a finite subcover. Finally, we take finitely
many barycentric subdivisions of $\p Y$ until the mesh is smaller
than the minimum diameter of the subcover.

Using the new triangulation, we repeat the assignment of the chain
complexes in each of the previous steps $0$ through $k$. In
particular, this defines $W^{k+1}_{\ell}(v)$ for all points in
$e(n-2-k)_{\ell}$ and by the modification argument indicated above,
we obtain a $(n-1)$-chain denoted by $W_{\ell}^{k+1}\times
e(n-2-k)_{\ell}$ for each $\ell\in E(n-2-k)$.

Once we have completed all the $n-1$ steps, we put together the
singular complexes constructed from the $(n-1)$-chains in each step
by means of their attaching maps along their common boundary which
will be explicitly computed. In order to finish the proof, we show
that $\p W^{n-1} = \p X_0 + \p X_1$ and thus once we attach $X_0$
and $X_1$ to the boundary we will have the cycle $\Gn$ we seek.

We consider the decomposition of $W^{n-1}$ from the $(n-1)$-chains
in each step, that is,
 \begin{equation}\label{eq:definition of cylinder as sum of complexes}
W^{n-1} = \sum_{k=0}^{n-2} ~\sum_{\ell\in E(n-2-k)}
W^{k+1}_{\ell}\times e(n-2-k)_{\ell}.
\end{equation}

For simplicity, let $S_k = \ds\sum_{\ell\in E(n-2-k)}
W^{k+1}_{\ell}\times e(n-2-k)_{\ell}$, then $\p W^{n-1}=
\ds\sum_{k=0}^{n-2} \p S_k$. We now analyze each term separately.
For $k=0$;
 \begin{equation}
\begin{split} \p S_{0} &= \left(\p\ds\sum_{\ell\in
E(n-2)} W^{1}_{\ell}\times e(n-2)_{\ell}\right)\\
&\!\!\!\!\!\!\!\!  = \sum_{\ell\in E(n-2)} \p W^{1}_{\ell}\times
e(n-2)_{\ell} + \ds\sum_{\ell\in E(n-2)} W^{1}_{\ell}\times \p
e(n-2)_{\ell}\\
&\!\!\!\!\!\!\!\!  =\sum_{\ell\in E(n-2)} (e(0)_{\ell}^0 +
e(0)_{\ell}^1)\times e(n-2)_{\ell} +\!\!\!\!\!\!\! \ds\sum_{\ell\in
E(n-2)}
\ds\sum_{j_{n-3}} W^{1}_{\ell}\times e(n-3)_{j_{n-3}}\\
&\!\!\!\!\!\!\!\!   = \quad\p X_0 + \p X_1 \quad + \quad
\ds\sum_{\ell\in E(n-2)}\ds\sum_{j_{n-3}} W^{1}_{\ell}\times
e(n-3)_{j_{n-3}},
\end{split}\nonumber
\end{equation}

\noindent where $j_{n-3}$ denotes the index of all $(n-3)$-cells
that belong to the boundary of $e(n-2)_{\ell}$. In general, for
$0<k<n-2$;
 \begin{equation}
\begin{split}
\p S_k &= \p\left(\ds\sum_{\ell\in E(n-2-k)} W^{k+1}_{\ell}\times
e(n-2-k)_{\ell}\right)\\
 &\!\!\!\!\!\!\!\! =\!\!\!\!\!\!\!\sum_{\ell\in E(n-2-k)} \p W^{k+1}_{\ell}\times
e(n-2-k)_{\ell} + \!\!\!\!\!\!\!\sum_{\ell\in E(n-2-k)}
W^{k+1}_{\ell}\times \p
e(n-2-k)_{\ell}\\
 &\!\!\!\!\!\!\!\! = \!\!\!\!\!\!\!\sum_{\ell\in E(n-2-k)}\sum_{j_{n-1-k}} W^{k}_{j_{n-1-k}}\times
e(n-1-k)_{\ell} \quad+\\
 &\!\!\!\!\!\!\!\! \quad\quad\quad\quad\quad+  \sum_{\ell\in
E(n-2-k)}\sum_{j_{n-3-k}}W^{k+1}_{\ell}\times e(n-3-k)_{j_{n-3-k}},
\end{split}\nonumber
\end{equation}

\noindent where $j_{n-1-k}$ denotes the index of all $(n-1-k)$-cells
that contain $e(n-2-k)_{\ell}$ and $j_{n-3-k}$ denotes the index of
all $({n-3-k})$-cells that belong to the boundary of
$e({n-3-k})_{\ell}$. Finally, for $k=n-2$;
 \begin{equation}
\begin{split}
\p S_{n-2}&= \p\left(\ds\sum_{\ell\in E(0)}
W^{n-1}_{\ell}\times e(0)_{\ell}\right)\\
 &=\sum_{\ell\in E(0)}\p W^{n-1}_{\ell}\times e(0)_{\ell}\quad +
\sum_{\ell\in E(0)}W^{n-1}_{\ell}\times \p e(0)_{\ell}\\
 &=\sum_{\ell\in
E(0)}\ds\sum_{j_1} W^{n-2}_{j_1}\times e(0)_{\ell}\quad + \quad0,
\end{split}\nonumber
\end{equation}

\noindent where $j_1$ denotes the index of all 1-cells that contain
$e(0)_{\ell}$ in its boundary.

Observe that the second summation term of $\p S_{k}$ is the same as
the first summation term of $\ S_{k+1}$ as we count each chain
twice. Since we are using $\Z_2$-coefficients, \eqref{eq:definition
of cylinder as sum of complexes} simplifies to $\p W^{n-1} = \p X_0
+ \p X_1$.

Finally, from the construction of $\W$ we see that for each
$q\in\text{supp\,}\W$, $q\in \finv(\Hv)$ for some $v\in\p Y$. This
concludes the proof of Lemma~\ref{lemma:Chain complex sausage of
revolution}. \QED


\section{Surjectivity}\label{sec:Proof of SUR}
In this section we consider the question of when a local
diffeomorphism $f:M\to\Rn$ is surjective, based on the topology of
the pre-images of hyperplanes. The trivial example of an inclusion
map of the region between two planes satisfies
Theorem~\ref{thm:injetividade general result for a manifold} but it
is not surjective. This indicates that further assumptions must be
added. On the other hand, we are able to eliminate the homological
assumption on the domain.


\begin{thm}\label{thm:surjectivity general result for a manifold}
Let $f:M\to\Rn$ be a local diffeomorphism where $M$ is a connected
manifold. If the pre-image of every affine hyperplane is non-empty
and acyclic, then $f$ is surjective.
\end{thm}

The proof is based on geometric constructions involving foliation
theory and the computation of linking numbers between certain
singular chain complexes in the range $\Rn$. We remark that since
the computation of linking numbers will occur in $\Rn$, it is not
necessary to make any further assumptions on the homology groups of
$M$. This is unlike the situation in Theorem~\ref{thm:injetividade
general result for a manifold}, where we assumed $H_{n-1}(M)=0$.

Combining Theorem~\ref{thm:surjectivity general result for a
manifold} and Theorem~\ref{thm:injetividade general result for a
manifold}, we obtain the following characterization of $\Rn$, for
$n\ge2$.

\begin{thm}\label{cor:Caracterizacao de Rn}
A smooth connected manifold $M$ is diffeomorphic to $\Rn$ if and
only if $H_{n-1}(M)=0$ and there exists a local diffeomorphism $f:M
\rightarrow \mathbb{R}^n$ such that the pre-image of every affine
hyperplane is non-empty and acyclic.
\end{thm}

\subsection{Proof of Theorem~\ref{thm:surjectivity general result for
a manifold}}

We establish surjectivity by showing that for each $R>0$ the ball of
radius $R$ is fully contained in $f(M)$, that is,
$\overline{B}(0;R)\subset f(M)$. Since translations do not change
any of our hypotheses, let us assume that $0\in f(M)$ and single out
$o\in\finv(0)\subset M$.

Next, from the local diffeomorphism assumption, there exists
$\eps>0$ such that $\overline{B}(0;\eps)\subset f(M)$ and
$\finv(B(0;\eps))$ has a diffeomorphic component $\We$ which
contains $o\in M$. Observe that for $R\leq\eps$, we trivially have
$\overline{B}(0;R)\subset f(M)$, so we restrict ourselves to the
case $R>\eps$.

We argue that we can find a way to expand $B(0;\eps)$ within the
image of $f$ so that it will contain a ball of radius any $R$. To
this end, we shall choose directions for this expansion as follows.
For $v\in S^{n-1}$, let $\Hv$ be the canonical codimension one
foliation of $\Rn$ by hyperplanes orthogonal to $v$. Since the leaf
space of $\Hv$ is homeomorphic to $\R$, we parameterize the leaves
of $\Hv$ by $\Hv(t)$ where $t$ is the distance of the hyperplane
$\Hv(t)$ to the origin. Because $\Hv(t) = \mathcal{H}_{-v}(-t)$, we
will only consider $t\ge0$. Let $\Nv=f^{\ast}\Hv$ be the pullback
foliation of $M$ which, by definition, has the connected components
of $\finv\left(\Hv(t)\right)$ as leaves. Since our hypotheses states
that the pre-images of hyperplanes are non-empty and connected, the
leaf space of $\Nv$ is homeomorphic to $\R$ and we then write
$\Nv(t)=\finv\left(\Hv(t)\right)$. Next, we claim that for each
$v\in S^{n-1}$, we may find a global transversal $\gv$ to the
foliation $\Nv$ that may be used to expand the image of $f$. More
precisely, we have the following result.

\begin{lemma}\label{lemma:global transversal in surjectivity} For
each $u\in \We$ with $f(u)=\eps v$, $v\in S^{n-1}$, there exists a
smooth path $\gv:[0,\infty) \rightarrow M$ with $\gv\pitchfork \Nv$
such that $f\big(\gv(t)\big)\in\Hv(t)$ for $t\in[0,R]$.
\end{lemma}

The proof follows directly from transverse modification arguments in
foliation theory (see \cite{Conlon}) and the fact that the leaves of
$\Nv$ are non-empty and connected.

As we proceed with the proof of Theorem~\ref{thm:surjectivity
general result for a manifold}, let us fix a canonical
identification of $S^{n-1}$ to $\p B(0;\eps)$ and $\We$. By applying
Lemma~\ref{lemma:global transversal in surjectivity}, we obtain
directions $\gv$ from which to expand $\We$ up to $\gv\cap\Nv(R)$.
Then for a fixed $v$, we can locally modify $\gv$ so that we can
carry an entire neighborhood of $v$ in $\We$ along $\gv$ using the
compactness of $\gv\bigl([0,R]\bigr)$. Repeating this process for
each $v\in\We$, we obtain a cover of $\We$. However, because there
is no canonical choice for $\gv$, points belonging to the
intersection of two neighborhoods may have multiply defined paths.
Our approach will be similar to the one in section \ref{sec:Proof of
INJ}. The key difference here is that we also need to control how
each neighborhood is pushed along the global transversal $\gv$.

Intuitively, we will push the cells of $\We$ along $\gv$ and
possibly create \textit{broken} pieces at each instant $t$. Using
the hypotheses that $\Nv(t)$ is acyclic and $f$ is a diffeomorphism
we will define bounding chains \emph{filling the gaps} in each leaf.
Furthermore, we will argue that these chain complexes constructed
for $s,t$ will be homologous, hence we say that they are homologous
relative to $t$. This process is depicted in Fig.~\ref{fig:intro to
expanding sphere} and it is stated precisely in the lemma below.

\begin{figure}[h]
\begin{center}
\ifthenelse{\equal{\pic}{1}}{
\begin{psfrags}
\psfrag{0}{$0$} \psfrag{Wn1eps}{$W^{n-1}(\eps)$}
\psfrag{Wn1R}{$W^{n-1}(R)$} \psfrag{gv}{$\gamma_v$}
\psfrag{Wn1t}{$W^{n-1}(t)$} \psfrag{p}{$q$}
\includegraphics{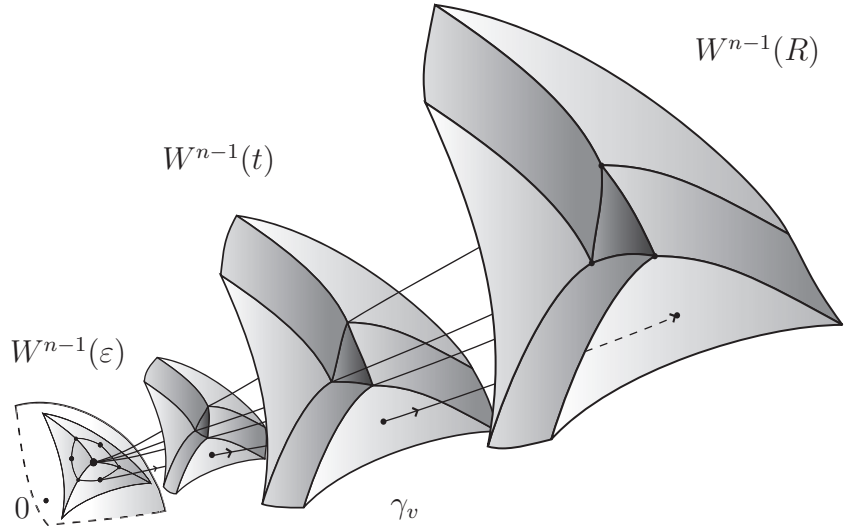}
\end{psfrags}
} {\[\fbox{\rule[2.5in]{2.5in}{0in}outline expanding spheres}\] }
\caption{A local assemblage of chain complexes based on the
triangulation of a sphere.} \label{fig:intro to expanding sphere}
\end{center}
\end{figure}

\begin{lemma}{}\label{lemma:spider web} For each $R>0$, there exists
a family of geometric (singular) chain complexes $\Wt$ that are
homologous in $M\backslash\set{o}$ for $t\in(0,R]$ such that:

\noi \emph{\textbf{i)}} For $t\in (0,\eps]$, $\Wt = \finv\bigl(\p
B(0;t)\bigr)$.

\noi \emph{\textbf{ii)}} If $q\in\text{supp\,}\Wt$, then $q\in
\Nv(t)$ for some $v\in S^{n-1}$.



\end{lemma}

The proof of Lemma~\ref{lemma:spider web} uses combinatorial
topology and foliation theory to explicitly construct such cycles.
Since the process is lengthy and rather technical, we postpone it
and continue with the proof of Theorem~\ref{thm:surjectivity general
result for a manifold}.

Our strategy to show that $\overline{B}(0;R)\subset f(M)$ is by
contradiction. Suppose there is $p\notin f(M)$. We compute the
linking number between $p$ and $f\big(\Wt\big)= \Zn(t)$ in two ways,
yielding different values. This argument is similar to standard
reasoning in degree theory and is geometric in nature. As before, we
work with $\Z_2$-coefficients.

Notice that since $f$ is continuous, we have that $\Zn(t)$ is a
family of homologous cycles in $\Rn\backslash\set{0}$. Then for
$t\in(0,\eps]$, $\Zn(t) = \p B(0;t)$ and we have that the origin is
contained in the \textit{inside} of $\Zn(\eps)$, more precisely, the
linking number between the origin and $\Zn(\eps)$ is equal to~1.


From (ii) of Lemma~\ref{lemma:spider web}, we have that
$0\notin\Zn(t)$ for each $t\in (0,R]$ and as mentioned above,
$\Zt\sim\Zr$ in $\Rn\backslash\set{0}$. Therefore as intersection
numbers, thus linking numbers are invariant under the same homology
class, we have
 \begin{equation}\label{eq:linking of the origin and Zn}
\vf(\Zn(t),0)=1\textrm{ for each }\eps<t\leq R,
\end{equation}

\noi where we consider the $0$-normal cycle formed by the origin and
a suitable point in the complement of a compact set
containing~$\Zt$.


We claim that $p$ is \textit{inside} $\Zr$, that is, $\vf(\Zr,p)=1$.
Indeed, consider the segment $Y^1$ from $0$ to $p$. We have that
$Y^1\cap \Zr = \emptyset$, otherwise it would imply $p\in f(M)$. By
definition,
 \begin{equation}\label{eq:disjoint bodies in proof of sobre thm}
\#(\Zn(R)\times Y^1) = \#(Y^1\times \Zn(R)) = 0. 
\end{equation}

Computing the linking number between the cycle $\Zr$ and $\p Y^1$
using the fact that $\p Y^1$ is a normal $0$-cycle, we have,
 \begin{equation}\label{eq:link of p with Zn}
 \begin{aligned}
0 & = \#(Y^1 \times \Zn(R)) = \vf(\p Y^1,\Zn(R)) = \vf(\Zn(R), \p Y^1)\\
  & = \vf(\Zn(R), 0-p) = \vf(\Zn(R), 0) - \vf(\Zn(R), p).
\end{aligned}
\end{equation}

Combining \eqref{eq:link of p with Zn} and $\vf(\Zn(R), 0)=1$ we
obtain
 \begin{equation}\label{eq:link of p with Zn equal to 1}
\vf(\Zn(R), p) =1.
\end{equation}

Observe that $\Ze=\p B(0;\eps)$, hence $\vf(\Ze,p) =0$. Finally,
from the assumption that $p\notin f(M)$, we have $\Zr\sim\Ze$ in
$\Rn-\{p\}$ and again by the invariance of linking numbers on the
homology class we obtain,
 \begin{equation}\label{eq:contradiction link}
\vf(\Zn(R), p) = \vf(\Zn(\eps), p)=0.
\end{equation}

This is a contradiction, therefore it must be the case that $p\in
f(M)$ and hence $f$ is surjective. \QED

\subsection{The Construction of a Family of Homologous
Cycles}\label{sec:The Construction of a spider web}

We now complete the proof of Theorem~\ref{thm:surjectivity general
result for a manifold} by establishing the technical proof of
Lemma~\ref{lemma:spider web}. We employ a similar technique as in
section~\ref{sec:Proof of INJ}, that is, we use triangulations as a
tool to keep track of intersections in the coverings. For
simplicity, let $\We =W$ and consider a triangulation of $W$ with
cells $e(k)_{\ell}$; $k=0,\ldots,n-1$ and $\ell\in E(k)\subset
\mathbb{N}$ where $E(k)$ is the set of indexes of all the $k$-cells
in $W$.


The idea of the construction of $\Wt$ is in essence geometric, and
can be outlined as follows. Consider a triangulation with
sufficiently small mesh. For each top dimensional cell of $W$, we
push it along a global transversal $\gv$ emanating from one of its
points up to the level $R$ and use the local modification of $\gv$
for points within the cell to push these points. The key issue is
that a point $v$ belonging to the boundary of a top dimensional cell
may be pushed along multiple choices of $\gv$, one for each top
dimensional cell it belongs to. Hence the $\Wt$ may not be well
defined; geometrically, this will create broken pieces at each
level. However, by considering the collections of cells that contain
$v$, via the link of $v$, we will show that for each $t$, the
multiply defined chain complexes form a cycle in the pre-image of
$\Hv(t)$. Thus by acyclicity, we can fill these gaps with bounding
chains. Furthermore, the process will be done so it is homologous
relative to $t$, that is, as we consider different chain complexes
for each $t$. Now, as we begin to formally describe $\Wt$, we will
do so in steps enumerated from $0$ to $n-1$, outlined below.

\jl

\noi\underline{Step 0:} For each point $v\in W$, suppose
$v\in\el{n-1}$ for some $\ell\in E(n-1)$. From
Lemma~\ref{lemma:global transversal in surjectivity} we obtain a
global transversal $\gv$ to the foliation $\Nv$. We then define the
following  0-chains, that is, points where $\gv$ intersect the
leaves of $\Nv$; Let $W^0_{\ell}(v,t)=\gv(t) \in \Nv(t)$ for $t\in
[\eps,R]$. By compactness of $\gv\bigl([0,R]\bigr)$ and the fact
that $f$ is a local diffeomorphism, there is a neighborhood
$\mathcal{V}_v\subset M$ of $\gv\bigl([0,R]\bigr)$ such that we can
continuously modify $\gv$ to obtain a global transversal
$\gamma_{v'}$ for all $v'$ in a neighborhood $\Ov\subset W$ of $v$,
as depicted in Fig~\ref{fig:Local modification of gv}.

\begin{figure}[h]
\begin{center}
\ifthenelse{\equal{\pic}{1}}{
\begin{psfrags}
\psfrag{gv}{$\gamma_v$} \psfrag{W0t}{$W^0_{\ell}(v,t)$}
\psfrag{W0R}{$W^0_{\ell}(v,R)$} \psfrag{v}{$v$}
\psfrag{Ov}{$\mathcal{O}_v$} \psfrag{0}{$o$}
\psfrag{Vv}{$\mathcal{V}_v$}
\includegraphics{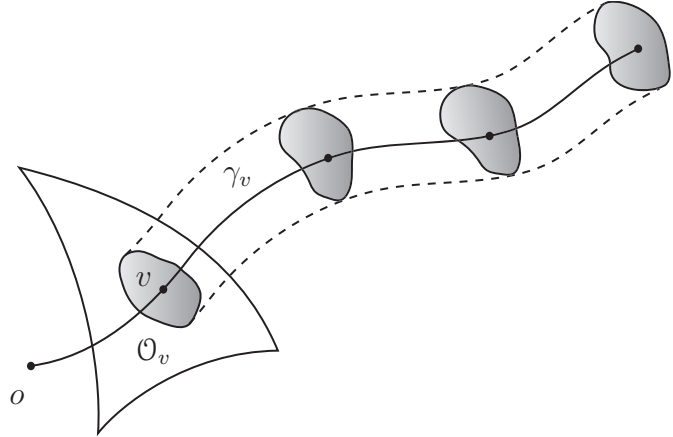}
\end{psfrags}
} {\[\fbox{\rule[0.75in]{3.0in}{0in}step0.eps}\] } \caption{Local
modification of global transversals.}\label{fig:Local modification
of gv}
\end{center}
\end{figure}

Then, by the compactness of $W$, we obtain a finite subcover of $W$
from $\set{\mathcal{O}_v}$. Now with such subcover, we iterate
finitely many barycentric subdivisions of $W$ until its mesh is
smaller than the minimum diameter of the subcover. In this process,
we obtain a new triangulation of $W$ with the property that, for
each $v\in \el{n-1}$, we may continuously define $\gv$ for all
points in $\el{n-1}$. In fact, we define, for each $t\in[\eps,R]$
and $\ell\in E(n-1)$, a $(n-1)$-chain denoted by
$W^0_{\ell}(t)\times\el{n-1}$ which is topologically equivalent to
$W^0_{\ell}(v,t)\times\el{n-1}$ and varies continuously on $t$.

\noi\underline{Step 1:} Consider points $v$ in the $(n-2)$-skeleton
of $W$. Suppose $v\in \el{n-2}$ for some $\ell\in E(n-2)$. Then $v$
belongs to the intersection of two $(n-1)$-cells that can be
determined by looking at the vertices in $\lk\left(\el{n-2}\right)$.
Without loss of generality let $v\in e(n-1)_1\cap e(n-1)_2$. Then
for each $t\in[\eps,R]$, we have defined in the previous step the
points $W^0_1(v,t), W^0_2(v,t)\in \Nv(t)$ along path emanating from
each top dimensional cell. Now we can join such points by a path
$W^1_{\ell}(v,t)$ lying in $\Nv(t)$ as it is acyclic.

Once we construct $W^1_{\ell}(v,t)$, we claim that it can be locally
modified for all points in a neighborhood of $(v,t)$ in
$W_{(n-2)}\times[\eps,R]$. Indeed, for each $v\in W_{(n-2)}$ and
$t\in [\eps, R]$, the path $W^1_{\ell}(v,t)$ is compact and hence we
may find a neighborhood $\U(v,t)\subset M$ of $W^1_{\ell}(v,t)$ such
that the (local) gradient flow of the height function $f_v:M\to\Rn$
given by $f_v(x) = \langle f(x),v \rangle$ can be used to
continuously define $W^1_{\ell}(v,t)$ for nearby $t$. Also, the fact
that $f$ is a local diffeomorphism continuously defines
$W^1_{\ell}(v,t)$ for all nearby $v$ in $W_{(n-2)}$.

The process above provides a cover $\left\{\U(v,t)\right\}$ of
$W_{(n-2)}\times [\eps,R]$. By compactness, we may find a finite
subcover which induces a cover of $W$. Indeed, in step 0 each top
dimensional cell is pushed diffeomorphically along the global
transversal $\gv$. We can now iterate finitely many barycentric
subdivisions of $W$ so its mesh is smaller than the minimum diameter
of the subcover above restricted to $W$. Next, we repeat all the
constructions up to this point using the new triangulation. Observer
that this guarantees that each cell $\el{n-2}$ is contained in a
member of the finite subcover.

Now let us consider a partition of $[\eps,R]$ induced by this
subcover, that is, we have $\eps = t_0<t_1<\cdots<t_N=R$ for some
$N\in \mathbb{N}$. From the choice of subdivision we can
continuously modify $W^1_{\ell}(v,t)$ for $t\in(t_i,t_{i+1})$ and
$v\in\el{n-2}$ by the argument above. The key problem is that for
the endpoint we may possibly have two chains defined, each coming
from the adjacent intervals. However, the fact that each leaf of
$\Nv$ is acyclic yields a similar construction as the transverse
modification method (see \cite{Conlon}) to ensure that whenever
ambiguity occurs, the choice will be homologous relative to $t$.

The details are as follows; For each $i=i,\ldots,N-1$, let the two
choices for a bounding chain $W^1_{\ell}(v,t_i)$ be
$W^1_{\ell}(v,t_i)^-$ and $W^1_{\ell}(v,t_i)^+$, where
$W^1_{\ell}(v,t_i)^-$ is the chain defined continuously from
$W^1_{\ell}(v,t)$, $t\in(t_{i-1}, t_i)$ and the second one,
$W^1_{\ell}(v,t_i)^+$, is the chain defined continuously from
$W^1_{\ell}(v,t)$, $t\in(t_{i}, t_{i+1})$. By default, we agree to
always choose $W^1_{\ell}(v,t_i)^+$. This will not be ambiguous
because we can choose either chain complex. Indeed,
$W^1_{\ell}(v,t_i)^-\sim\nolinebreak W^1_{\ell}(v,t_i)^+$ since from
construction they have the same boundary and, by acyclicity of
$\Nv(t_i)$, there is a bounding chain contained in $\Nv(t_i)$.
Finally, we must consider a new neighborhood $\widetilde{\U}(v,t_i)$
of such bounding chain where the restriction of $f$ is a
diffeomorphism. Doing so for every $v\in W_{(n-2)}$ and
$i=1,\ldots,N$ we obtain a new cover of $W$ by adding the collection
of sets $\widetilde{\U}(v,t_i)$ to the finite subcover considered up
to this point. This is done to ensure that the process will always
yield chains homologous relative to $t$. Then iterate finitely many
barycentric subdivisions of $W$ to obtain a triangulation with mesh
sufficiently small to define chains $W^1_{\ell}(v,t)$ continuously
for all points $v\in\el{n-2}$ for each $\ell$ and by construction
these chains are homologous relative to $t$ with continuously
varying bounding chains.

This adaptation of the transverse modification argument produces for
each $t\in[\eps,R]$ and $\ell\in E(n-2)$, a $(n-1)$-chain complex
denoted by $W^1_{\ell}(t)\times\el{n-2}$ which is topologically
equivalent to $W^1_{\ell}(v,t)\times\el{n-2}$, $v\in\el{n-2}$, and
is homologous relative to $t$. This concludes step 1.

Now we give the general procedure for $1<k\leq n-1$.

\noi\underline{Step $k$:} Consider points in the $(n-1-k)$-skeleton
of $W$. Suppose $v\in \el{n-1-k}$ for some $\ell\in E(n-1-k)$. We
are interested in identifying all the $(n-k)$-cells that contain
$\el{n-1-k}$ in its boundary, i.e., that contain $v$. This can be
accomplished by looking at the vertices of $\lk(\el{n-1-k})$. for
simplicity, suppose that those are $e(n-k)_1,
e(n-k)_2,\ldots,e(n-k)_m$ for some $m\in \mathbb{N}$. For each
$t\in[\eps,R]$, consider the $(k-1)$-chain;
 \begin{equation}\label{eq:boundary chain in step k of surjectivity}
W^{k-1}_1(v,t) + \cdots + W^{k-1}_m(v,t),
\end{equation}

\noi where the chains $W^{k-1}_j(v,t)$ were constructed in Step
$k-1$. We claim that the $(k-1)$-chain in~\eqref{eq:boundary chain
in step k of surjectivity} is a cycle. Indeed,
 \begin{equation}\label{eq:closedness of boundary chain in step k of surjectivity}
 \p\left( \sum_{j=1}^m W^{k-1}_j (v,t)\right) =
 \sum_{j=1}^m  \p W^{k-1}_j (v,t) =  \sum_{j=1}^m \left(
 \sum_{\ell'} W^{k-2}_{\ell'} (v,t)\right),
\end{equation}\nonumber

\noindent where $\ell'$ corresponds to the index of all
$(n-k+1)$-cells in $W$ that contain $\el{n-1-k}$. Now the argument
is completely analogous to the one given in the injectivity case,
i.e, it follows from the observation that an edge contains exactly
two vertices. Since $\Nv(t)$ is acyclic, there exists a $k$-chain
$W^{k}_{\ell}(v,t)$ that bounds $\ds\sum_{j=1}^m W^{k-1}_j (v,t)$.

We now argue that the chain $W^{k}_{\ell}(v,t)$ may be continuously
modified with respect to $v$ and within the same homology class
relative to $t$. This is similar to the construction as in step 1,
except that we repeat it for each interval of the partition obtained
in step $k-1$, so we omit the details. Finally, for each
$t\in[\eps,R]$ and $\ell\in E(n-1-k)$, we obtain a $(n-1)$-chain
denoted by $W^k_{\ell}(t)\times\el{n-1-k}$ which for
$v\in\el{n-1-k}$, is topologically equivalent to
$W^k_{\ell}(v,t)\times\el{n-1-k}$, and is homologous relative to the
parameter~$t$.

\begin{remark} In the last step $n-1$, we consider the $0$-skeleton of $W$ which
is discrete, hence no further subdivisions are necessary.
\end{remark}

Once all $n$ steps are completed, for each $t\in[\eps,R]$, we put
together all the constructed chain complexes via the obvious
attaching maps based on the intersection of the cells in $W$ as
indicated by their construction in each step. This defines $\Wt$ as
follows;
 \begin{equation}
W^{n-1}(t) =  \ds\sum_{k=0}^{n-1} \sum_{\ell\in E(n-1-k)}
W^{k}_{\ell}(t)\times e(n-1-k)_{\ell}. \nonumber
\end{equation}

We observe that $\Wt\sim\W(s)$ in $M\backslash\set{o}$ for
$t,s\in[\eps,R]$. Indeed, by construction if $t,s\in (t_{i-1},t_i)$
for $i=1,\ldots, N$, we can use the local gradient flow of the
corresponding height functions to continuously modify the chain
$\Wt$ to $\W(s)$, each chain $W^{n-1-k}_{\ell}(t)\times e(k)_{\ell}$
at a time. Otherwise, the only problem is at the end points $t_i$,
where again by construction, the chains of $\W(t_i)$ are homologous
in $\mathcal{N}(t_i)= \bigl\lbrace p\in M | p\in \Nv(t_i)\text{ for
some } v\in S^{n-1}\bigr\rbrace$, thus ensuring that $\Wt$ are
homologous in $M\backslash\set{o}$.

It now remains to show that $\Wt$ is a cycle in $M$. The computation
below is quite similar to the one done in section~\ref{sec:Proof of
INJ}. For simplicity, let $S_k(t) = \ds\sum_{\ell\in E(n-1-k)}
W^{k}_{\ell}(t)\times e(n-1-k)_{\ell}$, so that $\p W^{n-1}(t)
=\nolinebreak\ds \sum_{k=0}^{n-1} \p S_k(t)$.

Considering each term separately, we have for $k=0$;
 \begin{equation}
\begin{split}
\p S_{0}(t) &= \p\left(\ds\sum_{\ell\in E(n-1)}
W^{0}_{\ell}(t)\times e(n-1)_{\ell}\right) \\
& = \sum_{\ell\in E(n-1)} \p W^{0}_{\ell}(t)\times
e(n-1)_{\ell}\quad + \!\!\!\!\!\!\!\sum_{\ell\in E(n-1)}
W^{0}_{\ell}(t)\times \p
e(n-1)_{\ell}\\
& =\quad 0\quad+\quad\sum_{\ell\in E(n-1)}\ds\sum_{j_{n-2}}
W^{0}_{\ell}\times e(n-2)_{j_{n-2}},
\end{split}\nonumber
\end{equation}

\noindent  where $j_{n-2}$ denotes the index of all $(n-2)$-cells
that belong to the boundary of $e(n-1)_{\ell}$. In general, for
$1<k<n-1$;
 \begin{equation}
\begin{split}
\p S_k(t)& = \p\left(\ds\sum_{\ell\in E(n-1-k)}
W^{k}_{\ell}(t)\times e(n-1-k)_{\ell}\right)\\
 & =\sum_{\ell\in E(n-1-k)} \p W^{k}_{\ell}(t)\times
e(n-1-k)_{\ell}\quad +\\
 & \quad\quad\quad\quad\quad\quad+ \sum_{\ell\in E(n-1-k)}
W^{k}_{\ell}(t)\times\p
e(n-1-k)_{\ell}\\
 & = \sum_{\ell\in E(n-1-k)}\ds\sum_{j_{n-k}}
W^{k-1}_{j_{n-k}}(t)\times e(n-1-k)_{\ell}\quad+\\
 & \quad\quad\quad\quad\quad\quad+   \sum_{\ell\in
E(n-1-k)} \ds\sum_{j_{n-2-k}} W^{k}_{\ell}(t)\times
e(n-2-k)_{j_{n-2-k}},
\end{split}\nonumber
\end{equation}

\noindent where $j_{n-k}$ denotes the index of all $({n-k})$-cells
that contain $e(n-1-k)_{\ell}$ and $j_{n-2-k}$ denotes the index of
all $({n-2-k})$-cells that belong to the boundary of
$e(n-1-k)_{\ell}$. Finally, for $k=n-1$,
 \begin{equation}
\begin{split}
\p S_{n-1}(t) & = \p\left(\ds\sum_{\ell\in E(0)}
W^{n-1}_{\ell}(t)\times e(0)_{\ell}\right)\\
& = \sum_{\ell\in E(0)}\p W^{n-1}_{\ell}(t)\times e(0)_{\ell}\quad
+\quad \sum_{\ell\in E(0)} W^{n-1}_{\ell}(t)\times\p e(0)_{\ell}\\
& =\sum_{\ell\in E(0)}\ds\sum_{j_1} W^{n-2}_{j_1}(t)\times
e(0)_{\ell} \quad + \quad\quad0,
\end{split}\nonumber
\end{equation}

\noindent where $j_1$ denotes the index of all 1-cells that contain
$e(0)_{\ell}$ in its boundary.

As we sum the terms in $\W(t)$, we see that the second summation
term appearing in $\p S_k(t)$ coincided with the first summation in
$\p S_{k+1}(t)$. Indeed, these terms count the same objects twice,
and hence we obtain zero (recall that we are working with $Z_2$
coefficients). Therefore $W^{n-1}(t)$ is a cycle and this finishes
the proof of Lemma~\ref{lemma:spider web}.~\QED

\section{Final Remarks}\label{sec:Remarks}
Having obtained independent results on injectivity and surjectivity,
we note that our main result follows from
Theorem~\ref{thm:injetividade general result for a manifold} and
Theorem~\ref{thm:surjectivity general result for a manifold}.
Recall,

\vskip0.1in

\noindent \textbf{Theorem~\ref{thm:characterization result}.} {\em A
local diffeomorphism $f:\mathbb{R}^n\to\mathbb{R}^n$ is bijective if
and only if the pre-image of every affine hyperplane is non-empty
and acyclic.}

\vskip0.1in

We also have the following analytic condition that establishes
whether a local diffeomorphism is bijective. Given a complete
Riemannian metric $g$ on $\Rn$ and a smooth function $h:\Rn\to\R$,
the gradient of $h$ relative to $g$, denoted by $\n^{g}h$, satisfies
$g_x(\n^{g}h,w) = d h_x(w)$ for all $w\in\Rn$. Our analytic result
is the following.

\begin{corollary} A local
diffeomorphism $f:\Rn\to\Rn$ is bijective if for each $v\in
S^{n-1}$, there exists a complete metric $g_v$ on $\Rn$ such that,
 \begin{equation}\label{eq:condition of analytic Corollary}
\inf_{x\in\Rn} |\nabla^{g_v} f_v(x)|_{g_v} >0.
\end{equation}
\end{corollary}

It is easy to see that such condition implies that the pre-images of
hyperplanes are acyclic, hence the result follows. We compare this
result with the work in \cite{NolletXavier1} where we can now choose
the metric to suit the unit vector. Finally, We can also state an
analytical result implying only injectivity.

\begin{corollary}\label{thm:Result with diff metrics for injectivity} A local
diffeomorphism $f:\mathbb{R}^n\rightarrow \mathbb{R}^n$ is injective
provided there exists $w\in S^{n-1}$ with the property that for each
unit vector $v$ perpendicular to $w$, there exists a complete
Riemannian metric $g_v$ on $\mathbb{R}^n$ such that,
 \begin{equation}\label{eq:vector field in diff metrics for injectivity}
 \inf_{x\in\mathbb{R}^n} \|Df(x)^{*}v\|_{g_v}>0.
\end{equation}

\end{corollary}

Observe that in our injectivity results we did not need topological
hypotheses on the pre-image of every hyperplanes, hence the result
holds.


\end{document}